\documentclass[twoside]{amsart}

\usepackage{amsmath,amssymb,amsfonts,amsthm,latexsym}
\usepackage[margin=1.1in]{geometry}
\usepackage{graphicx}
\usepackage[dvipsnames]{xcolor}
\usepackage{enumitem}
\usepackage{subcaption}
\usepackage{hyperref}
\usepackage{ytableau}
\usepackage{mathdots}
\usepackage{tikz-cd}
\usepackage{orcidlink}
\usepackage{mathtools}

\newtheorem{theorem}{Theorem}[section]
\newtheorem{lemma}[theorem]{Lemma}
\newtheorem{proposition}[theorem]{Proposition}
\newtheorem{corollary}[theorem]{Corollary}

\theoremstyle{definition}

\newtheorem{example}[theorem]{Example}
\theoremstyle{remark}
\newtheorem{remark}[theorem]{Remark}

\numberwithin{equation}{section}

\newcommand{\C}{\mathbb{C}}

\newcommand{\R}{\mathbb{R}}
\newcommand{\T}{\mathbb{T}}

\newcommand{\Z}{\mathbb{Z}}

\newcommand{\U}{\operatorname{U}}

\newcommand{\Rank}{\operatorname{rank}}
\newcommand{\diag}{\operatorname{diag}}

\newcommand{\bs}{\boldsymbol}

\newcommand{\ab}{\text{ab}}
\newcommand{\cl}{\text{cl}}
\newcommand{\OO}{\operatorname{O}}

\newcommand{\GL}{\operatorname{GL}}

\newcommand{\HOM}{\operatorname{Hom}}

\DeclarePairedDelimiter\abs{\lvert}{\rvert}%
\DeclarePairedDelimiter\norm{\lVert}{\rVert}%
\makeatletter
\let\oldabs\abs
\def\abs{\@ifstar{\oldabs}{\oldabs*}}
\let\oldnorm\norm
\def\norm{\@ifstar{\oldnorm}{\oldnorm*}}
\makeatother

\begin{document}

\title{Higher orbifold Euler characteristics of torus representations}

\author[N.~Peters]{Nathaniel Peters}
\address{Department of Mathematics and Statistics,
Skidmore College,
815 North Broadway,
Saratoga Springs, NY 12866, USA}

\author[C.~Seaton]{Christopher Seaton\,\orcidlink{0000-0001-8545-2599}}
\address{Department of Mathematics and Statistics,
Skidmore College,
815 North Broadway,
Saratoga Springs, NY 12866, USA}
\email{cseaton@skidmore.edu}


\thanks{Both authors were supported by the Skidmore College Faculty Student Summer Research Program.
C.S. was supported by an AMS-Simons Research Enhancement Grant for PUI Faculty.}

\begin{abstract}
We compute the $\Gamma$-Euler characteristic of a finite-dimensional representation of the compact torus $\mathbb{T}^k$
where $\Gamma$ is an arbitrary finitely presented group. We additionally compute the $\Gamma$-Euler characteristic
of the unit sphere in such a representation as well as the corresponding symplectic quotient at a nonzero level of the moment map.
In each case, the $\Gamma$-Euler characteristic is expressed in terms of the weight matrix for the representation as well
as a description of the abelianization of $\Gamma$.
\end{abstract}

% xxxxxxxxxxxxxxxxxxxxxxxxxxxxxxxxxxxxxxxxxxxxxxxxxxxxxxxxxxxxxxxxxxxxxxxxx
% xxxxxxxxxxxxxxxxxxxxxxxxxxxxxxxxxxxxxxxxxxxxxxxxxxxxxxxxxxxxxxxxxxxxxxxxx
% xxxxxxxxxxxxxxxxxxxxxxxxxxxxxxxxxxxxxxxxxxxxxxxxxxxxxxxxxxxxxxxxxxxxxxxxx

\setcounter{tocdepth}{3}
\maketitle

\tableofcontents

% xxxxxxxxxxxxxxxxxxxxxxxxxxxxxxxxxxxxxxxxxxxxxxxxxxxxxxxxxxxxxxxxxxxxxxxxx
% xxxxxxxxxxxxxxxxxxxxxxxxxxxxxxxxxxxxxxxxxxxxxxxxxxxxxxxxxxxxxxxxxxxxxxxxx
% xxxxxxxxxxxxxxxxxxxxxxxxxxxxxxxxxxxxxxxxxxxxxxxxxxxxxxxxxxxxxxxxxxxxxxxxx

\section{Introduction}
\label{sec:intro}

The \emph{$\Gamma$-Euler characteristics} $\chi_\Gamma$ are topological invariants generalizing the various orbifold Euler characteristics
\cite{AtiyahSegal,DixonHarveyEA,FarsiSeaGenOrbEuler,HirzebruchHoefer,Roan,TamanoiMorava,TamanoiCovering} to groupoids that need not have finite isotropy groups.
Each finitely presented group $\Gamma$ yields such an invariant, which is defined for \emph{orbit space definable groupoids}, topological groupoids whose
isotropy groups are compact Lie and whose orbit spaces are relatively tame, reflected by being contained in an o-minimal structure.
The $\Gamma$-Euler characteristics were introduced in the case of compact Lie group actions and $\Gamma = \Z^\ell$ in \cite{G-ZLM-Hcompact},
called \emph{orbifold Euler characteristic of order $\ell$}, and in the general context in \cite{FarSeaJLMS}.

In \cite{FarsiMeitSeaton}, the authors began a program to study these invariants and their dependence on the groupoid structure
by computing them for large classes of example. Focusing on the case of a translation groupoid $G\ltimes V$ where $G$ is a compact Lie group
and $V$ is a finite-dimensional $G$-representation, \cite{FarsiMeitSeaton} includes computations of $\chi_\Gamma(G\ltimes V)$ where $G$ is the circle or real
orthogonal group $\OO(2, \R)$. In addition, the $\Gamma$-Euler characteristic is computed for the symplectic quotient associated to a unitary representation
of an arbitrary compact Lie group $G$.

In this paper, we continue this program by computing the $\Gamma$-Euler characteristics of $\T^k$-spaces where $\T^k$ is the compact $k$-dimensional torus.
Specifically, we consider the case of a finite-dimensional representation of $\T^k$, the unit sphere in such a representation, and the symplectic quotient
of a unitary representation at a nonzero level of the moment map. Using the fact that $\T^k$ is abelian, we are able to determine explicit formulas
for the $\Gamma$-Euler characteristic for arbitrary finitely presented $\Gamma$, each expressed in terms of the weight matrix $A$ of the action as
well as a description of the abelianization $\Gamma^{\ab}$ of $\Gamma$; see Theorem~\ref{thrm:SmithAbelianization}.
For symplectic quotients, many but not all of the cases addressed by Theorem~\ref{thrm:SymplecticGamEC} are orbifolds,
in which case the $\Gamma$-Euler characteristic coincides with the $\Gamma$-orbifold Euler characteristic of
\cite{FarsiSeaGenOrbEuler,TamanoiMorava,TamanoiCovering}. The cases addressed by Theorems~\ref{thrm:RepGamEC} and \ref{thrm:SphereGamEC} are
only orbifolds in trivial examples.

The outline of this paper is as follows. In Section~\ref{sec:Back}, we recall relevant background and fix notation. Section~\ref{sec:GamECFulLRep}
computes $\chi_\Gamma$ for the case of a representation; see Theorem~\ref{thrm:RepGamEC} and Remark~\ref{rem:RepRealFactor}.
Section~\ref{sec:GamECSubspaces} considers subspaces of such a representation and computes $\chi_\Gamma$ of the unit sphere in
Section~\ref{subsec:GamECSubspacesSpheres} and of a symplectic quotient in
Section~\ref{subsec:GamECSubspacesSymplectic}; see Theorems~\ref{thrm:SphereGamEC} and \ref{thrm:SymplecticGamEC}.

% xxxxxxxxxxxxxxxxxxxxxxxxxxxxxxxxxxxxxxxxxxxxxxxxxxxxxxxxxxxxxxxxxxxxxxxxx
% xxxxxxxxxxxxxxxxxxxxxxxxxxxxxxxxxxxxxxxxxxxxxxxxxxxxxxxxxxxxxxxxxxxxxxxxx
% xxxxxxxxxxxxxxxxxxxxxxxxxxxxxxxxxxxxxxxxxxxxxxxxxxxxxxxxxxxxxxxxxxxxxxxxx

\section*{Acknowledgements}

This paper developed from N.P.'s research project conducted through the Skidmore College Faculty Student Summer Research Program
and Skidmore Mathematics and Statistics Department, and the authors gratefully acknowledge the support of the department and college for these activities.
Both authors were supported by the Skidmore College Faculty Student Summer Research Program, and
C.S. was supported by an AMS-Simons Research Enhancement Grant for PUI Faculty.
C.S. would like to thank Carla Farsi and Emily Proctor for helpful discussions.

% xxxxxxxxxxxxxxxxxxxxxxxxxxxxxxxxxxxxxxxxxxxxxxxxxxxxxxxxxxxxxxxxxxxxxxxxx
% xxxxxxxxxxxxxxxxxxxxxxxxxxxxxxxxxxxxxxxxxxxxxxxxxxxxxxxxxxxxxxxxxxxxxxxxx
% xxxxxxxxxxxxxxxxxxxxxxxxxxxxxxxxxxxxxxxxxxxxxxxxxxxxxxxxxxxxxxxxxxxxxxxxx

\section{Background}
\label{sec:Back}

Here we collect the relevant background material and establish the notation we will use throughout this paper.
We review the representation theory of the real torus in Section~\ref{subsec:BackTorusReps} and refer the reader to
\cite[Ch.~II, Sec.~8]{BroeckerTomDieck} for more details. Section~\ref{subsec:BackGammaEC} recalls the definition of the
$\Gamma$-Euler characteristic and related background, see \cite{FarSeaJLMS} for more details, as well as auxiliary results from \cite{FarsiMeitSeaton}.
In Section~\ref{subsec:BackSmith}, we recall the Smith normal form of an integer matrix \cite{SmithLinCongr} as well as two
applications that will play a large role in the paper. For more details about Smith normal form over an arbitrary PID, we
refer the reader to the modern treatment of \cite[Sec.~5.3]{AdkinsWeintraub}. The second application, which computes the abelianization
of a finitely presented group, is detailed in \cite[Sec.~9.2]{HoltEickOBrienBook}.

Throughout this paper, it will be helpful to adopt the conventions that $0^0 = 1$ and that the determinant of a non-square matrix,
including a $k\times 0$ matrix with $k > 1$, is zero.
The latter convention is merely a useful book-keeping device that will allow us to write several equations uniformly.

% xxxxxxxxxxxxxxxxxxxxxxxxxxxxxxxxxxxxxxxxxxxxxxxxxxxxxxxxxxxxxxxxxxxxxxxxx

\subsection{Representations of Tori}
\label{subsec:BackTorusReps}

Let $k \geq 1$ and let $\T^k$ denote the $k$-dimensional real torus, the product of $k$ copies of the circle $\T^1 = \{ z\in\C : \abs{z} = 1\}$.
By \cite[Ch.~II, Prop.~8.5]{BroeckerTomDieck}, every finite-dimensional real representation $U$ of $\T^k$ can be expressed as the
sum of the underlying real representation of a unitary representation $V\simeq\C^n$ and a trivial real representation $W\simeq\R^d$.
Moreover, there is a basis for $V$ with respect to which the $\T^k$-action is diagonalizable and can be described by an
integer \emph{weight matrix} $A = (a_{i,j}) \in\Z^{k\times n}$. Specifically, there are coordinates $(\theta_1,\ldots,\theta_k)$
for the Lie algebra $\mathfrak{t}^k\simeq\R^k$ of $\T^k$ and complex coordinates $(v_1,\ldots,v_n)$ for $V$ such that the action
of $\exp(\bs{\theta}) = \exp(\theta_1,\ldots,\theta_k)\in\T^k$ on $\mathbf{v} = (v_1,\ldots,v_n)\in V$ is given by
$\exp(\bs{\theta}A)\mathbf{v}$, where $\bs{\theta}$ is treated as a $1\times k$ row matrix and $\exp(\bs{\theta}A)\in\GL(V,\C)\simeq\GL(n,\C)$.
Note that there is a Hermitian structure on $V$ with respect to which $\exp(\bs{\theta}A)\in\U(V)\simeq\U(n,\C)$;
however, we focus on the underlying real structure of $V$, and use complex coordinates only for convenience.

Row-reducing $A$ with row operations invertible over $\Z$, i.e., elements of $\GL(k,\Z)$, corresponds to changing the basis of $\mathfrak{t}^k$ and hence
the product decomposition of $\T^k$ but does not change the representation. Hence, we will sometimes assume that $A$ is
upper-triangular. This can always be accomplished via integer row operations by putting $A$ in \emph{Hermite normal form};
see \cite[Sec.~5.2]{AdkinsWeintraub} for more information.
Permuting the columns of $A$ merely permutes the coordinates of $V$. Though we will not use this fact, we note that
multiplying a column by $-1$ changes the unitary representation but does not change the underlying real representation.

We let $\rho\colon \T^k \to \GL(V,\C)$ denote the representation homomorphism $\rho(\exp(\bs{\theta})) = \exp(\bs{\theta}A)$.
If $T\simeq\T^n$ denotes the maximal compact torus in $\GL(V,\C)$ consisting of diagonal matrices whose diagonal entries have modulus $1$, then
the image $\rho(\T^k)$ of $\T^k$ in $\GL(V,\C)$ is a subtorus of $T$ of dimension $\Rank A$. Note that we do not assume that $k\leq n$
nor that $A$ has full rank; that is, we do not assume that $V$ is faithful.

Let $U = V\oplus W$ be a finite-dimensional $\T^k$-representation as above.
For much of this paper, the real factor $W$ will play a trivial role, so we will focus on the unitary representation $V$.
In Section~\ref{subsec:GamECSubspacesSpheres}, where the real factor deserves more attention, we will reduce to the cases that $\dim_\R W = 0$ or $1$ as follows.
If $\dim_\R W$ is even, then $W$ can be written as a direct sum of copies of $\R^2$ which, via a choice of basis, can be given a Hermitian structure.
Hence, we can redefine $V$ so that $U\simeq V$ as real representations; this amounts to adding $(\dim_\R W)/2$ zero columns to the matrix $A$.
If $\dim_\R W$ is odd, then we do the same to a real hyperspace of $W$ and are left with
$V\oplus W^\prime$ where $V$ is unitary and $\dim_\R W^\prime = 1$. We will use complex coordinates for $V\simeq\C^n$ as above and real coordinates
$(w)$ for $W\simeq\R$.

The $\T^k$-isotropy group of a point $\mathbf{v}\in V$ is determined by its zero coordinates, which we describe with the following notation.
Let $[n] = \{ 1, 2, \ldots, n\}$. For $I\subseteq[n]$, we let
\[
    V_I =   \{ (v_1,\ldots,v_n) \in V : v_i \neq 0 \Leftrightarrow i\in I \}.
\]
We let $\abs{I}$ denote the cardinality of $I$, let $A_I$ denote the $(k\times\abs{I})$-minor of $A$ consisting of the columns of $A$ indexed by elements
of $I$, and let $r_I = \Rank A_I$. Then the closure of $V_I$ is given by
\[
    V_I^{\cl} =   \{ (v_1,\ldots,v_n) \in V : v_i \neq 0 \Rightarrow i\in I \}.
\]
Note that $V_I^{\cl}$ is a unitary $\T^k$-representation with weight matrix $A_I$. We let $\rho_I\colon \T^k\to\GL(V_I,\C)$,
and then the image of $\rho_I$ is a torus of dimension $r_I$.
Finally, we let $G_I$ denote the common $\T^k$-isotropy group of points in $V_I$, equivalently, the kernel of $\rho_I$.

% xxxxxxxxxxxxxxxxxxxxxxxxxxxxxxxxxxxxxxxxxxxxxxxxxxxxxxxxxxxxxxxxxxxxxxxxx

\subsection{The $\Gamma$-Euler Characteristic}
\label{subsec:BackGammaEC}

An \emph{o-minimal structure} on $\R$ is a sequence of distinguished subsets of $\R^n$ for each $n$. An example is the collection of \emph{semialgebraic sets},
those given by finite unions of sets defined by finitely many polynomial equations and inequalities. The reader unfamiliar with o-minimal structures
can read this paper considering only the o-minimal structure of semialgebraic sets. See \cite{vandenDriesBook} for more details and background on o-minimal structures.

Fix an o-minimal structure on $\R$ that contains the semialgebraic sets. We say that an element of this o-minimal structure is a \emph{definable set},
and a function from a definable set $A\subseteq\R^m$ to a definable set $B\subseteq\R^n$ is \emph{definable} if its graph is a definable subset of $\R^{m + n}$.
An \emph{affine definable space} is a topological space $X$ equipped with a topological embedding $\iota_X\colon X\to \R^n$ whose image is definable,
and a subset $C\subseteq X$ is a \emph{definable subset} if $\iota_X(C)$ is definable. A \emph{morphism of affine definable spaces} is a continuous function $f$
from an affine definable spaces $(X,\iota_X)$ to affine definable spaces $(Y, \iota_Y)$ such that $\iota_Y\circ f\circ\iota_X^{-1}$ is definable. 

An affine definable space $X$ admits a finite definable cell decomposition and triangulation, potentially including open cells or simplices,
and hence has a well-defined Euler characteristic $\chi(X)$ given by the alternating sum of the number of cells or simplices of each
dimension; see \cite[Ch.~3, 4, 8]{vandenDriesBook}. By \cite[Thrm.~2.2]{BekeInvarEC}, the Euler characteristic of $X$ depends only on the underlying topological space
of $X$. Note that $\chi$ is finitely additive and multiplicative but is not homotopy invariant.

The $\Gamma$-Euler characteristics were defined in \cite{FarSeaJLMS} for \emph{orbit space definable groupoids}, a class of groupoids whose orbit spaces are affine
definable spaces. In this paper, we will only consider translation groupoids $G\ltimes X$ given by the action of a compact Lie group $G$ on an affine definable space.
Hence, for the benefit of the reader unfamiliar with groupoids, we restrict to this context as we recall the definition of the $\Gamma$-Euler characteristics.
See also \cite[Def.~2.1 and 2.2]{G-ZLM-Hcompact} for a formulation of $\chi_{\Z^\ell}$ of a $G$-space, which corresponds to the
\emph{orbifold Euler characteristic of order $\ell$} in this reference, and \cite[Sec.~4]{FarsiProctorSeatonBurnside} for the relation to the \emph{universal
Euler characteristic}.

Let $G$ be a compact Lie group. By \cite[Ch.~VI, Sec.~XII, Thrm.~4]{ChevalleyThryLieGrps}, $G$ admits a finite-dimensional representation, from which it inherits the
structure of an affine definable space with respect to which the product and inverse operations are definable.
If $X$ is an affine definable space and as well a $G$-space such that the action map $G\times X\to X$ is definable, we say that $X$ is an \emph{ affine definable $G$-space}.
If $U$ is a real finite-dimensional linear representation of $G$, then $U$ is an affine definable $G$-space, as is any $G$-stable definable subset of $U$. As well,
the corresponding homomorphism $G\to\GL(U,\R)$ is a morphism of affine definable spaces.
By \cite[Ch.~10, Cor.~2.18]{vandenDriesBook}, as $G$ is compact, $X$ admits a \emph{definably proper quotient}; this means that the orbit space $G\backslash X$
is an affine definable space, the quotient map $X\to G\backslash X$ is a morphism of affine definable spaces, and the quotient map is as well \emph{definably proper},
i.e., the preimage of a definable compact set is compact. See \cite[Sec.~10.2]{vandenDriesBook} and \cite[Sec.~2.2]{FarSeaJLMS} for more details.

For $x\in X$, let $G_x\leq G$ denote the isotropy group of $x$ and let $Gx\subseteq X$ denote the $G$-orbit of $x$.
We assume that the partition of $G\backslash X$ into points whose isotropy groups are isomorphic as topological groups is a finite partition
of $G\backslash X$ into definable sets. If $\Gamma$ is a finitely presented group, then the \emph{$\Gamma$-Euler characteristic of $G\ltimes X$}
is defined to be
\begin{equation}
\label{eq:DefGamEC}
    \chi_\Gamma(G\ltimes X)
    =
    \int_{G\backslash X} \chi\big(G_x\backslash\HOM(\Gamma,G_x)\big) \, d\chi(G x),
\end{equation}
where the $G_x$-action on $\HOM(\Gamma,G_x)$ is by conjugation and the integral is with respect to the Euler characteristic;
see \cite{ViroIntegralEulerChar}.
For more details or the more general formulation of $\chi_\Gamma$, see \cite[Def.~3.8]{FarSeaJLMS}.

We recall the following from \cite[Eq.~(2.3)]{FarsiMeitSeaton}.
If $G\backslash X = Y_1\sqcup\cdots\sqcup Y_s$ is a partition of $G\backslash X$ into definable sets such that for each $i$, the isotropy groups of
points with orbits in $Y_i$ are isomorphic as topological groups to a fixed topological group $G_i$, then
\begin{equation}
\label{eq:GamECSum}
    \chi_\Gamma(G\ltimes X)
    =
    \sum\limits_{i=1}^s \chi(Y_i)\chi\big(G_i\backslash\HOM(\Gamma, G_i)\big).
\end{equation}
In the cases under consideration, the groups $G_i$ will always be abelian so that the conjugation action on $\HOM(\Gamma, G_i)$ is trivial;
in this case Equation~\eqref{eq:GamECSum} becomes
\begin{equation}
\label{eq:GamECSumAbel}
    \chi_\Gamma(G\ltimes X)
    =
    \sum\limits_{i=1}^s \chi(Y_i)\chi\big(\HOM(\Gamma, G_i)\big).
\end{equation}

Finally, we recall the following, which is a formulation of \cite[Thrm.~3.1]{FarsiMeitSeaton} in the specific context where it will be applied here.

\begin{proposition}
\label{prop:TorusActionSubgroupEC}
Let $T$ be a torus and let $X$ be an affine definable $T$-space. If $H$ is a closed subgroup of $T$,
then $\chi(H\backslash X) = \chi\big( (H\backslash X)^{T/H}\big)$ where $(H\backslash X)^{T/H}$ denotes the set of points
fixed by $T/H$.
\end{proposition}
\begin{proof}
As $T$ is abelian, the action of $T$ on $X$ induces an action of $T/H$ on $H\backslash X$. As this action is the
composition of the definable action map $T\times X\to X$ with the definable quotient map $X\to H\backslash X$, it is definable. The claim then follows from
\cite[Thrm.~3.1]{FarsiMeitSeaton} and the fact that $T/H$ is a torus.
\end{proof}

We will most frequently use this proposition in the case that $\dim H < \dim T$ and $X^T = \emptyset$,
implying that $\chi\big( (H\backslash X)^{T/H}\big) = \chi(\emptyset) = 0$.

% xxxxxxxxxxxxxxxxxxxxxxxxxxxxxxxxxxxxxxxxxxxxxxxxxxxxxxxxxxxxxxxxxxxxxxxxx

\subsection{Smith Normal Form of Integer Matrices}
\label{subsec:BackSmith}

We now recall the Smith normal form of an integer matrix as well as two applications that will be important in the sequel.

Let $A\in\Z^{k\times n}$ be an integer matrix. For a positive integer $i$,
if $\Rank A \leq i$, we let $d_i(A)$ denote the $\gcd$ of the set of $i\times i$ minors of $A$.
If $i > \Rank A$, we set $d_i(A) = 0$, and we let $d_0(A) = 1$. The integer $d_i(A)$ is called the \emph{$i$\textsuperscript{th} determinantal divisor of $A$};
see \cite[Ch.5, Def.~3.10]{AdkinsWeintraub}.

\begin{theorem}[Smith Normal Form {\cite{SmithLinCongr},\cite[Ch.~5, Thrm.~3.1 and Lem.~3.11]{AdkinsWeintraub}}]
\label{thrm:SmithNormForm}
Let $A\in\Z^{k\times n}$ of rank $r_A$. There are matrices $P\in\GL(k,\Z)$ and $Q\in\GL(n,\Z)$ such that
\begin{equation}
\label{eq:SmithNormForm}
    PAQ =   \begin{pmatrix} D & \bs{0}_{r_A\times(n-r_A)} \\ \bs{0}_{(k-r_A)\times r_A} & \bs{0}_{(k-r_A)\times(n-r_A)} \end{pmatrix} \in\Z^{k\times n},
\end{equation}
where $D = \diag(s_1,\ldots,s_{r_A})\in\Z^{r_A\times r_A}$ is a diagonal $r_A\times r_A$ matrix with each $s_i\neq 0$,
$s_i\mid s_{i+1}$ for $1\leq i\leq r_A - 1$, and for $i = 1,\ldots, r_A$,
\[
    s_i =   \frac{d_i(A)}{d_{i-1}(A)}.
\]
\end{theorem}

The block matrix on the right side of Equation~\eqref{eq:SmithNormForm} is called the \emph{Smith normal form of $A$}.
The entries $s_1,\ldots,s_{r_A}$ are the \emph{invariant factors of $A$}; they are uniquely determined up to sign, see
\cite[Ch.~5, Rem.~3.2(2)]{AdkinsWeintraub}, and we will adopt the convention that they are always taken to be positive.
Hence, the zero matrix has no invariant factors.

An important application of Theorem~\ref{thrm:SmithNormForm} for us will be the following.
Let $A\in\Z^{k\times n}$ be the weight matrix of an $n$-dimensional unitary $\T^k$-representation $V\simeq\C^n$ and choose corresponding
coordinates $(v_1,\ldots,v_n)$ for $V$. By the \emph{kernel} of the representation, we mean the kernel of $\rho$, equivalently,
the subgroup of $\T^k$ that acts trivially on $V$ and the isotropy group of a point $(v_1,\ldots,v_n)\in V$ such that each $v_i$ is nonzero.

\begin{proposition}
\label{prop:SmithIsotropy}
Let $A\in\Z^{k\times n}$ be the weight matrix of an $n$-dimensional unitary $\T^k$-representation $V$. The kernel $K\leq\T^k$ of the representation $V$
is isomorphic to $\Z/s_1\Z \oplus\cdots\oplus \Z/s_{r_A}\Z\oplus\T^{k-r_A}$ where $r_A = \Rank A$ and the $s_i$ are the invariant factors of $A$.
\end{proposition}
\begin{proof}
Let $S$ denote the Smith normal form of $A$ so that there are $P\in\GL(k,\Z)$ and $Q\in\GL(n,\Z)$ such that $PAQ = S$.
For an element of $\T^k$ of the form $\exp(\bs{\theta})$ with $\bs{\theta}\in\mathfrak{t}^k$,
the action on $V$ in coordinates corresponding to $A$ is given by left-multiplication by $\exp(\bs{\theta}A) = \exp(\bs{\theta}P^{-1}SQ^{-1})$.
The matrix $P^{-1}$ corresponds to a change of basis in $\mathfrak{t}^k$ that preserves the lattice $2\pi\Z^k$,
and the matrix $Q^{-1}$ corresponds to a change of basis in $\mathfrak{t}^n$ preserving the lattice $2\pi\Z^n$, i.e., both preserve the kernel of the
corresponding exponential map. Specifically, let $T\simeq\T^n$ denote the maximal compact torus in $\GL(n,\C)$ of diagonal matrices with complex units on the diagonal,
and note that $\exp(\bs{\theta}A)\in T\leq\GL(n,\C)$ is the identity if and only if
$\bs{\theta}A = \bs{\theta}P^{-1}SQ^{-1}\in 2\pi\Z^n$.
As $Q$ preserves $2\pi\Z^n$, this is equivalent to $\bs{\theta}P^{-1}S\in 2\pi\Z^n$. Then as $P$ preserves $2\pi\Z^k$, this is the case if and only if
$\bs{\theta}P^{-1}$ is in the kernel of the representation $\exp(\bs{\phi})\mapsto\exp(\bs{\phi}S)$, i.e., the representation with weight matrix $S$.
Using coordinates corresponding to the weight matrix $S$, an element $\exp(\bs{\phi})$ is in the kernel if and only if
$\exp(\bs{\phi}) = (\zeta_1,\ldots,\zeta_{r_A}, t_{r_A+1},\ldots, t_k)$ where each $\zeta_i$ is an $s_i$\textsuperscript{th} root of unity and the $t_i$
are arbitrary elements of $\T^1$.
\end{proof}

The reader is cautioned that the Proposition~\ref{prop:SmithIsotropy} cannot be used to determine the isotropy groups of individual points in $V$.
It is possible that points in the representation corresponding to $A$ have isotropy groups that do not occur in the representation corresponding to $S$.

\begin{example}
\label{ex:SmithUnequalIsotropy}
Let $A = \begin{pmatrix} 1 & 2 \end{pmatrix}$ be the weight matrix for a $\T^1$-action on $V = \C^2$, and then points of the form $(0,v_2)$ with $v_2\neq 0$
have isotropy $\{\pm 1\}$ isomorphic to $\Z/2\Z$. The Smith normal form of $A$ is given by the decomposition
\[
    A   =   \begin{pmatrix} 1 \end{pmatrix}
            \begin{pmatrix} 1 & 0 \end{pmatrix}
            \begin{pmatrix} 1 & 2 \\ 0 & 1 \end{pmatrix},
\]
and the only isotropy groups of points in the representation with weight matrix $S = \begin{pmatrix} 1 & 0 \end{pmatrix}$ are the trivial group and $\T^1$.
Proposition~\ref{prop:SmithIsotropy} tells us only that $A$ and $S$ correspond to representations with the same kernel; in this case, both kernels are trivial.

Of course, the $\Z/2\Z$-isotropy in $V$ can be detected by considering the trivial Smith normal form of the submatrix $A_{\{2\}} = \begin{pmatrix} 2 \end{pmatrix}$
of $A$.
\end{example}

Another application of Theorem~\ref{thrm:SmithNormForm} we will encounter is the following, which, because $\T^k$ is abelian,
will allow us to reduce the computation of the $\Gamma$-Euler characteristic to the case where $\Gamma$ is abelian.
If $\Gamma$ is a finitely generated group, we let $\Gamma^{\ab} = \Gamma/[\Gamma, \Gamma]$ denote the \emph{abelianization} of $\Gamma$, where
$[\Gamma, \Gamma]$ is the commutator subgroup. It is well known that if $H$ is an abelian group, then there is an isomorphism between the abelian groups
$\HOM(\Gamma^{\ab},H)$ and $\HOM(\Gamma,H)$ given by composing a group homomorphism $\phi\colon \Gamma^{\ab}\to H$ with the quotient homomorphism
$\Gamma\to\Gamma/[\Gamma,\Gamma] = \Gamma^{\ab}$. From the definition, it follows that this isomorphism preserves the image in $H$.
In \cite[Sec.~9.2]{HoltEickOBrienBook}, it is demonstrated how to use Smith normal form to compute the abelianization of a finitely presented group;
along with the above observations, this yields the following.

\begin{theorem}[Abelianization of $\Gamma$ {\cite[Prop.~2.68 and Sec.~9.2]{HoltEickOBrienBook}}]
\label{thrm:SmithAbelianization}
Let $\Gamma$ be a finitely presented group with presentation $\langle \gamma_1,\ldots,\gamma_N, R_1,\ldots,R_q\rangle$. For each relation $R_i$
of the form $\gamma_{i_1}^{p_{i_1}}\cdots\gamma_{i_u}^{p_{i_u}}$, let $R_i^+$ denote the column vector in $\Z^N$ whose $j$\textsuperscript{th}
entry is the sum $\sum_{i_\ell = j} p_{i_\ell}$. Let $M = (R_1^+ \cdots R_q^+)\in\Z^{N\times q}$
denote the matrix with columns $R_i^+$, let $r_M = \Rank M$, and let $\sigma_1,\ldots,\sigma_{r_M}$ denote the invariant factors of $M$. Then
\[
    \Gamma^{\ab}    \simeq
    \Z^{N - r_M} \oplus \bigoplus\limits_{i=1}^{r_M} \Z/\sigma_i\Z.
\]
If $H$ is an abelian group, then composition with the quotient homomorphism $\Gamma\to\Gamma/[\Gamma,\Gamma]$ induces an isomorphism of abelian groups
\[
    \HOM(\Gamma, H)
    \simeq
    \HOM(\Gamma^{\ab}, H)
    \simeq
    \HOM(\Z, H)^{N - r_M} \oplus \bigoplus\limits_{i=1}^{r_M} \HOM(\Z/\sigma_i\Z, H).
\]
\end{theorem}

% xxxxxxxxxxxxxxxxxxxxxxxxxxxxxxxxxxxxxxxxxxxxxxxxxxxxxxxxxxxxxxxxxxxxxxxxx
% xxxxxxxxxxxxxxxxxxxxxxxxxxxxxxxxxxxxxxxxxxxxxxxxxxxxxxxxxxxxxxxxxxxxxxxxx
% xxxxxxxxxxxxxxxxxxxxxxxxxxxxxxxxxxxxxxxxxxxxxxxxxxxxxxxxxxxxxxxxxxxxxxxxx

\section{$\Gamma$-Euler Characteristic of a Representation}
\label{sec:GamECFulLRep}

Let $V$ be a finite-dimensional unitary representation of $\T^k$. In this section, we compute the $\Gamma$-Euler characteristic of
$\T^k\ltimes V$ for an arbitrary finitely presented group $\Gamma$; see Theorem~\ref{thrm:RepGamEC}. A direct consequence is the
$\Gamma$-Euler characteristic of $\T^k\ltimes U$ where $U$ is a finite-dimensional real representation; see Remark~\ref{rem:RepRealFactor}.
These results generalize \cite[Thrm.~4.1 and Cor.~4.3]{FarsiMeitSeaton}, which address the case $k = 1$.

We begin with two lemmas computing the Euler characteristic of pieces of the orbit space of the form $\T^k\backslash V_I$.

\begin{lemma}
\label{lem:RepRank<|I|}
Let $A\in\Z^{k\times n}$ be the weight matrix of a unitary $\T^k$-representation $V$.
If $I\subseteq[n]$ such that $r_I  < \abs{I}$, then $\chi\big(\T^k\backslash V_I\big) = 0$.
\end{lemma}
\begin{proof}
Assume $r_I < \abs{I}$, which implies that $I\neq\emptyset$, and recall that $\rho_I\colon\T^k\to\GL(V_I^{\cl},\C)$ denotes the representation
with weight matrix $A_I$,
the restriction of the $\T^k$-representation on $V$ to $V_I$. Let $T\simeq\T^{\abs{I}}$ denote the maximal compact torus of $\GL(V_I^{\cl},\C)$
given by diagonal matrices whose diagonal elements have modulus $1$, and note that the action of $T$ on $V_I^{\cl}$ restricts to an
action on $V_I$. Then $\rho_I(\T^k)$ is a closed subgroup of $T$ isomorphic to $\T^{r_I}$ so that
$T/\rho_I(\T^k)$ is a positive-dimensional torus. As the $i$\textsuperscript{th} coordinate of each point in $V_I$ is nonzero,
$V_I^T = \emptyset$, and the claim follows from Proposition~\ref{prop:TorusActionSubgroupEC}.
\end{proof}

Note that the hypothesis $r_I  < \abs{I}$ of Lemma~\ref{lem:RepRank<|I|} always holds when $\abs{I} > k$.

\begin{lemma}
\label{lem:RepRank=|I|}
Let $A\in\Z^{k\times n}$ be the weight matrix of a unitary $\T^k$-representation $V$.
If $I\subseteq[n]$ such that $r_I  = \abs{I}$, then $\chi\big(\T^k\backslash V_I\big) = (-1)^{\abs{I}}$.
\end{lemma}
\begin{proof}
If $I = \emptyset$, then $r_I = 0 = \abs{I}$ and $V_I = \{\mathbf{0}\}$ so that $\T^k\backslash V_I$ is a point and the claim holds.
So assume $I\neq\emptyset$ such that $r_I = \abs{I}$.
By permuting coordinates in $V$, we assume for simplicity that
$I = [r_I]$ for $1\leq r_I \leq n$, and then $r_I = \abs{I} \leq k$. By changing coordinates in $\mathfrak{t}^k$, we assume that
$A_I$ is upper-triangular; e.g., we may put $A$ in Hermite normal form \cite[Sec.~5.2]{AdkinsWeintraub}.
Then the rank hypothesis implies that the diagonal entries of $A_I$ are nonzero. Let $(v_1,\ldots,v_{r_I},0,\ldots,0)\in V_I$. For $t_1\in\T^1$, we have
\[
    \rho_I(t_1,1,\ldots,1) (v_1,\ldots,v_{r_I},0,\ldots,0) = (t_1^{a_{1,1}}v_1, t_1^{a_{1,2}}v_2,\ldots,t_1^{a_{{r_I},1}}v_{r_I},0,\ldots,0).
\]
As $a_{1,1}\neq 0$ we may choose $t_1\in\T^1$ such that $t_1^{a_{1,1}}v_1$ lies on the positive real axis $\R_+$. For $t_2\in\T^1$, we have
\[
    \rho_I(1,t_2,1,\ldots,1)(t_1^{a_{1,1}}v_1, t_1^{a_{1,2}}v_2,\ldots,t_1^{a_{{r_I},1}}v_{r_I},0,\ldots,0)
    =   (t_1^{a_{1,1}}v_1, t_2^{a_{2,2}}t_1^{a_{1,2}}v_2,\ldots,t_2^{a_{2,{r_I}}}t_1^{a_{{r_I},1}}v_{r_I},0,\ldots,0),
\]
so as $a_{2,2}\neq 0$, we may choose $t_2$ such that $t_2^{a_{2,2}}t_1^{a_{1,2}}v_2\in\R_+$.
Continuing in this way, there is an element $\mathbf{t}\in\T^k$ such that $\rho_I(\mathbf{t})(v_1,\ldots,v_n)\in(\R_+)^{r_I}\times\{\mathbf{0}\}\subset V_I$.
As the $\T^k$-action preserves the modulus of each coordinate, distinct points in $(\R_+)^{r_I}\times\{\mathbf{0}\}\subset V_I$ are elements of distinct
$\T^k$-orbits so that the quotient map $V_I\to\T^k\backslash V_I$ restricts to a continuous bijection $(\R_+)^{r_I}\times\{\mathbf{0}\}\to\T^k\backslash V_I$.
As the quotient map is definable, this restriction is as well definable. The claim then follows from the fact that the Euler
characteristic is preserved under definable bijections; see \cite[Ch.~4, Prop.~2.4]{vandenDriesBook}.
\end{proof}

We now use Lemmas~\ref{lem:RepRank<|I|} and \ref{lem:RepRank=|I|} to reduce the computation of $\chi_\Gamma\big(\T^k\ltimes V\big)$ to the computation of the
$\chi\big(\HOM(\Gamma, G_I)\big)$. Recall that $G_I$ denotes the isotropy group of a point in $V_I$

\begin{lemma}
\label{lem:RepGenGamEC}
Let $A\in\Z^{k\times n}$ be the weight matrix of a unitary $\T^k$-representation $V$ and let $\Gamma$ be a finitely presented group. Then
\begin{equation}
\label{eq:RepGenGamEC}
    \chi_\Gamma\big(\T^k\ltimes V\big)
    =
    \sum\limits_{\substack{I\subseteq[n] \\ r_I = \abs{I}}} (-1)^{\abs{I}} \chi\big(\HOM(\Gamma, G_I)\big).
\end{equation}
\end{lemma}
\begin{proof}
Partition $V = \bigsqcup_{I\subseteq[n]} V_I$. Applying Equation~\eqref{eq:GamECSumAbel}, to this partition, we have
\[
    \chi_\Gamma\big(\T^k\ltimes V\big)
    =
    \sum\limits_{I\subseteq[n]} \chi\big(\T^k\backslash V_I\big)\chi\big(\HOM(\Gamma, G_I)\big).
\]
For each $I\subseteq[n]$, clearly, $r_I \leq \abs{I}$.
If $r_I < \abs{I}$, we have $\chi\big(\T^k\backslash V_I\big) = 0$ by Lemma~\ref{lem:RepRank<|I|} and the term vanishes.
If $r_I = \abs{I}$, then by Lemma~\ref{lem:RepRank=|I|}, $\chi\big(\T^k\backslash V_I\big) = (-1)^{\abs{I}}$.
\end{proof}

By Lemma~\ref{lem:RepGenGamEC}, it remains only to compute $\chi\big(\HOM(\Gamma, G_I)\big)$, which we now do with the help of
Proposition~\ref{prop:SmithIsotropy} and Theorem~\ref{thrm:SmithAbelianization}. Recall that we adopt the conventions that $0^0 = 1$
and the determinant of a non-square matrix, including a $k\times 0$ matrix, is $0$.

\begin{lemma}
\label{lem:RepGenGamHom}
Let $A\in\Z^{k\times n}$ be the weight matrix of a unitary $\T^k$-representation $V$ and let $\Gamma$ be a finitely presented group.
Choose a finite presentation of $\Gamma$ and let $N$ denote the number of generators and $r_M$ and $\sigma_1,\ldots,\sigma_{r_M}$ denote
the rank and invariant factors, respectively, of the matrix $M$ defined in Theorem~\ref{thrm:SmithAbelianization}.
For each $I\subseteq[n]$,
\begin{equation}
\label{eq:RepGenGamHom}
    \chi\big(\HOM(\Gamma, G_I)\big)
    =
    \abs{\det A_I}^{N - r_M}
    \prod\limits_{i=1}^{r_M} \sigma_i^{k - r_I} \prod\limits_{j=1}^{r_I} \gcd(\sigma_i, s_{I,j}),
\end{equation}
where $r_I = \Rank A_I$ and $s_{I,1}, \ldots, s_{I,r_I}$ are the invariant factors of $A_I$.
\end{lemma}
\begin{proof}
Fix $I\subseteq[n]$ and note that as $A$ has rank at most $k$, $r_I \leq k$.
As $G_I \leq \T^k$ is abelian, we have by Theorem~\ref{thrm:SmithAbelianization} that
\begin{equation}
\label{eq:RepGenGamHomDecomp1}
    \chi\big(\HOM(\Gamma, G_I)\big)
    =
    \chi\Big(\HOM(\Z, G_I)^{N - r_M} \oplus \bigoplus\limits_{i=1}^{r_M} \HOM(\Z/\sigma_i\Z, G_I)\Big).
\end{equation}
The map $\HOM(\Z, G_I)\to G_I$ that sends each $\phi\in\HOM(\Z, G_I)$ to $\phi(1)$ is an isomorphism, so by the multiplicativity of $\chi$,
Equation~\eqref{eq:RepGenGamHomDecomp1} becomes
\begin{equation}
\label{eq:RepGenGamHomDecomp2}
    \chi\big(\HOM(\Gamma, G_I)\big)
    =
    \chi(G_I)^{N - r_M} \prod\limits_{i=1}^{r_M} \chi\big(\HOM(\Z/\sigma_i\Z, G_I)\big).
\end{equation}

First assume $r_M < N$, i.e., the abelianization $\Gamma^{\ab}$ of $\Gamma$ has infinite order. If $r_I < k$,
then by Proposition~\ref{prop:SmithIsotropy}, $G_I$ is the product of a finite group and a positive dimensional torus.
Hence, again by the multiplicativity of $\chi$, the first factor vanishes and $\chi\big(\HOM(\Gamma, G_I)\big) = 0$.
As well, $A_I$ is either not square when $\abs{I} < k$ or square and singular when $\abs{I} = k$;
in either case, as $N - r_M > 0$, we have $\abs{\det A_I}^{N - r_M} = 0$.

So suppose $r_I = k$. Again by Proposition~\ref{prop:SmithIsotropy}, $G_I\simeq \bigoplus_{j=1}^{r_I} \Z/s_{I,j}\Z$ so that $G_I$ has order
$\prod_{j=1}^{r_I} s_{I,j}$. It follows that $\chi(G_I) = \prod_{j=1}^{r_I} s_j$, which by Theorem~\ref{thrm:SmithNormForm} is equal to
\begin{equation}
\label{eq:RepZECProductInvarFactors}
    \prod\limits_{j=1}^{r_I} d_j(A_I)/d_{j-1}(A_I)
    =   d_k(A_I)
    = \abs{\det A_I}.
\end{equation}
For each $i=1,\ldots,r_M$, we have isomorphisms
\[
    \HOM(\Z/\sigma_i\Z, G_I)
    \simeq
    \bigoplus\limits_{j=1}^{r_I}
    \HOM(\Z/\sigma_i\Z, \Z/s_{I,j}\Z),
\]
and each $\HOM(\Z/\sigma_i\Z, \Z/s_{I,j}\Z)$ contains $\gcd(\sigma_i, s_{I,j})$ elements, one sending the generator of $\Z/\sigma_i\Z$ to each of the
$\gcd(\sigma_i, s_{I,j})$ elements of order divisible by $\sigma_i$. It then follows from Equation~\eqref{eq:RepGenGamHomDecomp2} that
\begin{align*}
    \chi\big(\HOM(\Gamma, G_I)\big)
    &=      \abs{\det A_I}^{N - r_M} \prod\limits_{i=1}^{r_M} \prod\limits_{j=1}^{r_I} \gcd(\sigma_i, s_{I,j})
    \\&=    \abs{\det A_I}^{N - r_M} \prod\limits_{i=1}^{r_M} \sigma_i^{k - r_I} \prod\limits_{j=1}^{r_I} \gcd(\sigma_i, s_{I,j}),
\end{align*}
where $\sigma_i^{k - r_I} = 1$ as $\abs{I} = k$.

Now assume $r_M = N$ so that $\Gamma^{\ab} \simeq \bigoplus_{i=1}^{r_M} \Z/\sigma_i\Z$ is finite and the factor $\chi(G_I)^{N - r_M}$
on the right side of Equation~\eqref{eq:RepGenGamHomDecomp2} is $1$. Proposition~\ref{prop:SmithIsotropy} yields
$G_I\simeq \T^{k - r_I} \oplus\bigoplus_{j=1}^{r_I} \Z/s_{I,j}\Z$ so that for each $i=1,\ldots,r_M$,
\[
    \HOM(\Z/\sigma_i\Z, G_I)
    \simeq
    \HOM(\Z/\sigma_i\Z,  \T)^{k - r_I}\oplus \bigoplus\limits_{j=1}^{r_I} \HOM(\Z/\sigma_i\Z, \Z/s_{I,j}\Z).
\]
An element of $\HOM(\Z/\sigma_i\Z, \T)$ can send a generator of $\Z/\sigma_i\Z$ to any of the $\sigma_i$\textsuperscript{th} roots of unity
so that $\HOM(\Z/\sigma_i\Z, \T)$ contains $\sigma_i$ elements. As noted above, each $\HOM(\Z/\sigma_i\Z, \Z/s_{I,j}\Z)$ contains $\gcd(\sigma_i, s_{I,j})$
elements so that
\[
    \chi\big(\HOM(\Z/\sigma_i\Z, G_I)\big)
    =   \sigma_i^{k - r_I} \prod\limits_{j=1}^{r_I} \gcd(\sigma_i, s_{I,j}).
\]
As $r_M = N$ implies $\abs{\det A_I}^{N - r_M} = 1$,
\begin{align*}
    \chi\big(\HOM(\Gamma, G_I)\big)
    &=      \prod\limits_{i=1}^{r_M} \sigma_i^{k - r_I} \prod\limits_{j=1}^{r_I} \gcd(\sigma_i, s_{I,j})
    \\&=    \abs{\det A_I}^{N - r_M} \prod\limits_{i=1}^{r_M} \sigma_i^{k - r_I} \prod\limits_{j=1}^{r_I} \gcd(\sigma_i, s_{I,j}),
\end{align*}
completing the proof.
\end{proof}

Combining Lemmas~\ref{lem:RepGenGamEC} and \ref{lem:RepGenGamHom} yields the following, the main result of this section.

\begin{theorem}
\label{thrm:RepGamEC}
Let $A\in\Z^{k\times n}$ be the weight matrix of a unitary $\T^k$-representation $V$ and let $\Gamma$ be a finitely presented group.
Choose a finite presentation of $\Gamma$ and let $N$ denote the number of generators and $r_M$ and $\sigma_1,\ldots,\sigma_{r_M}$ denote
the rank and invariant factors, respectively, of the matrix $M$ defined in Theorem~\ref{thrm:SmithAbelianization}.
Then
\begin{equation}
\label{eq:RepGamEC}
    \chi_\Gamma\big(\T^k\ltimes V\big)
    =
    \sum\limits_{\substack{I\subseteq[n] \\ \abs{I} = r_I}} (-1)^{\abs{I}} \abs{\det A_I}^{N - r_M}
        \prod\limits_{i=1}^{r_M} \sigma_i^{k - \abs{I}} \prod\limits_{j=1}^{\abs{I}} \gcd(\sigma_i, s_{I,j}),
\end{equation}
where $r_I = \Rank A_I$ and $s_{I,1}, \ldots, s_{I,r_I}$ are the invariant factors of $A_I$.
\end{theorem}

\begin{remark}
\label{rem:RepGamECSimplifications}
Note that if $r_M < N$, i.e., $\Gamma^{\ab}$ is infinite, then $\abs{\det A_I}^{N - r_M} = 0$ unless $A_I$ is square, i.e., $\abs{I} = k$.
Then Equation~\eqref{eq:RepGamEC} simplifies to
\begin{equation}
\label{eq:RepGamECGamABInfinite}
    \chi_\Gamma\big(\T^k\ltimes V\big)
    =
    (-1)^k \sum\limits_{\substack{I\subseteq[n] \\ \abs{I} = k = r_I}} \abs{\det A_I}^{N - r_M}
        \prod\limits_{i=1}^{r_M} \prod\limits_{j=1}^{k} \gcd(\sigma_i, s_{I,j}).
\end{equation}
If $r_M = N$ so that $\Gamma^{\ab}$ is finite, then $\abs{\det A_I}^{N - r_M} = 1$ for each $I\subseteq[n]$ so that Equation~\eqref{eq:RepGamEC} becomes
\begin{equation}
\label{eq:RepGamECGamABFinite}
    \chi_\Gamma\big(\T^k\ltimes V\big)
    =
    \sum\limits_{\substack{I\subseteq[n] \\ \abs{I} = r_I}} (-1)^{\abs{I}}
        \prod\limits_{i=1}^{r_M} \sigma_i^{k - \abs{I}} \prod\limits_{j=1}^{\abs{I}} \gcd(\sigma_i, s_{I,j}).
\end{equation}
\end{remark}

\begin{remark}
\label{rem:RepRealFactor}
Theorem~\ref{thrm:RepGamEC} also yields the $\Gamma$-Euler characteristic of an arbitrary real finite-dimensional representation $U$ of $\T^k$ as follows.
By \cite[Prop.~8.5]{BroeckerTomDieck}, $U\simeq V\oplus W$ where $V$ is a unitary representation of $\T^k$ and $W$ is a trivial real representation of $\T^k$;
see Section~\ref{subsec:BackTorusReps}. Then $\T^k\ltimes U$ is isomorphic as a topological groupoid to the product $(\T^k\ltimes V)\times W$ where
$W$ is the base groupoid with only unit arrows, and $\chi_\Gamma(W) = \chi(W) = (-1)^{\dim W}$ for any finitely presented
$\Gamma$. By \cite[Lem.~4.17]{FarSeaJLMS}, it follows that $\chi_\Gamma\big(\T^k\ltimes U\big) = (-1)^{\dim W}\chi_\Gamma\big(\T^k\ltimes V\big)$.
\end{remark}

Let us note a few important cases of Theorem~\ref{thrm:RepGamEC}. If $\Gamma = \Z^\ell$, then using the standard presentation of $\Gamma$,
$N = \ell$ and the matrix $M$ in Theorem~\ref{thrm:SmithAbelianization} is the zero matrix. Hence $r_M = 0$ so that the product over $i$ in
Equation~\eqref{eq:RepGamEC} or \eqref{eq:RepGamECGamABInfinite} is empty, and we obtain the following.

\begin{corollary}
\label{cor:RepZ^ellEC}
Let $A\in\Z^{k\times n}$ be the weight matrix of a unitary $\T^k$-representation $V$. Then for any positive integer $\ell$,
\begin{equation}
\label{eq:RepZ^ellEC}
    \chi_{\Z^\ell}(\T^k\ltimes V)
    =
    (-1)^k\sum\limits_{\substack{I\subseteq[n] \\ \abs{I} = k}} \abs{\det A_I}^\ell.
\end{equation}
In particular,
\begin{equation}
\label{eq:RepZEC}
    \chi_\Z(\T^k\ltimes V)
    =
    (-1)^k\sum\limits_{\substack{I\subseteq[n] \\ \abs{I} = k}} \abs{\det A_I}.
\end{equation}
\end{corollary}

Similarly, if $\Gamma = \Z/m$ for a positive integer $m$, then $N = 1$ and the matrix $M$ in Theorem~\ref{thrm:SmithAbelianization} is $(m)$,
trivially in Smith normal form. Hence $r_M = 1$ and $\sigma_i = m$, yielding the following.

\begin{corollary}
\label{cor:RepZ/mEC}
Let $A\in\Z^{k\times n}$ be the weight matrix of a unitary $\T^k$-representation $V$. For any positive integer $m$,
\begin{equation}
\label{eq:RepZ/mEC}
    \chi_{\Z/m\Z}(\T^k\ltimes V)
    =
    \sum\limits_{\substack{I\subseteq[n] \\ \abs{I} = r_I}} (-1)^{\abs{I}} m^{k - \abs{I}}
        \prod\limits_{i=1}^{\abs{I}} \gcd(m, s_{I,i}),
\end{equation}
where $r_I = \Rank A_I$ and $s_{I,1}, \ldots, s_{I,r_I}$ are the invariant factors of $A_I$.
\end{corollary}

% xxxxxxxxxxxxxxxxxxxxxxxxxxxxxxxxxxxxxxxxxxxxxxxxxxxxxxxxxxxxxxxxxxxxxxxxx
% xxxxxxxxxxxxxxxxxxxxxxxxxxxxxxxxxxxxxxxxxxxxxxxxxxxxxxxxxxxxxxxxxxxxxxxxx
% xxxxxxxxxxxxxxxxxxxxxxxxxxxxxxxxxxxxxxxxxxxxxxxxxxxxxxxxxxxxxxxxxxxxxxxxx

\section{$\Gamma$-Euler Characteristic of a Subspace of a Representation}
\label{sec:GamECSubspaces}

We now use the techniques of the previous section to compute $\chi_\Gamma\big(\T^k\ltimes X\big)$ for certain $G$-stable definable subsets $X$ of a $\T^k$-representation.
In Section~\ref{subsec:GamECSubspacesSpheres}, we consider the case where $X$ is the unit sphere in an arbitrary real representation, while in
Section~\ref{subsec:GamECSubspacesSymplectic}, we consider the symplectic quotient, where $X$ is a fiber of the moment map in a unitary representation.

% xxxxxxxxxxxxxxxxxxxxxxxxxxxxxxxxxxxxxxxxxxxxxxxxxxxxxxxxxxxxxxxxxxxxxxxxx

\subsection{The Unit Sphere}
\label{subsec:GamECSubspacesSpheres}

Let $U$ be a finite-dimensional real $\T^k$-representation. As explained in Section~\ref{subsec:BackTorusReps}, we can assume without loss of generality
that there is a unitary representation $V\simeq\C^n$ such that as real representations, either $U\simeq V$ or $U\simeq V\oplus W$ where $\dim_\R W = 1$.
In the first case, we continue to use complex coordinates $(v_1,\ldots,v_n) \in\C^n$ for $V$, while in the second, we use coordinates
$(v_1,\ldots,v_n,w)$ where each $v_i\in\C$ and $w\in\R$.

To discuss both cases simultaneously, we let $d = \dim_\R U \bmod 2$ so that $\dim_\R U = d + \dim_\R V$.
Let $S^{2n + d - 1}$ denote the unit sphere in $U$, i.e., $S^{2n + d - 1}$ is given by $\{ (v_1,\ldots,v_n) : \sum_{i=1}^n \abs{v_i}^2 = 1\}$ when $d = 0$ and
$\{ (v_1,\ldots,v_n, w) : w^2 + \sum_{i=1}^n \abs{v_i}^2 = 1\}$ when $d = 1$.
For each $I\subseteq[n]$, let
\[
    S_{I}^{2n + d - 1} =
    \begin{cases}
        S^{2n + d - 1}\cap V_I,         &   d = 0,
        \\
        S^{2n + d - 1}\cap V_I\times W, &   d = 1.
    \end{cases}
\]
We first prove two lemmas that compute $\chi\big(\T^k\backslash S_I^{2n + d - 1}\big)$. The proofs are similar to those of
Lemmas~\ref{lem:RepRank<|I|} and \ref{lem:RepRank=|I|}.

\begin{lemma}
\label{lem:SphereRank<|I|}
Let $A\in\Z^{k\times n}$ be the weight matrix of a unitary $\T^k$-representation $V$ and
let $U$ be a finite-dimensional real $\T^k$-representation with $U\simeq V$ or $U\simeq V\oplus W$ as above, where $\dim_\R W = 1$.
If $I\subseteq[n]$ such that $r_I  < \abs{I}$, then $\chi\big(\T^k\backslash S_I^{2n + d - 1}\big) = 0$.
\end{lemma}
\begin{proof}
Assume $r_I < \abs{I}$ so that $I\neq\emptyset$.
Let $T\simeq\T^{\abs{I}}$ denote the torus in $\GL(V_I^{\cl},\C)$ given by diagonal matrices whose diagonal elements have modulus $1$.
As the action of $T$ on $V_I$ preserves the modulus of each coordinate in $V$, $T$ acts on $S_I^{2n + d - 1}$ when $d = 0$, as does
$T\times\{1\}\simeq T$ when $d = 1$.
As the coordinates of points in $S_I^{2n + d - 1}$ corresponding to $i\in I \neq\emptyset$ are nonzero, there are no points in
$S_I^{2n + d - 1}$ fixed by $T$ when $d = 0$ or $T\times\{1\}$ when $d = 1$.
Then as $\rho_I(\T^k)$ is a subgroup of $T$ isomorphic to a torus of dimension $r_I < \abs{I}$, an application of
Proposition~\ref{prop:TorusActionSubgroupEC} completes the proof.
\end{proof}

\begin{lemma}
\label{lem:SphereRank=|I|}
Let $A\in\Z^{k\times n}$ be the weight matrix of a unitary $\T^k$-representation $V$ and
let $U$ be a finite-dimensional real $\T^k$-representation with $U\simeq V$ or $U\simeq V\oplus W$ as above, where $\dim_\R W = 1$.
If $I\subseteq[n]$ is nonempty such that $r_I  = \abs{I}$, then $\chi\big(\T^k\backslash S_I^{2n + d - 1}\big) = (-1)^{r_I + d - 1}$.
If $I = \emptyset$ so that $r_I = \abs{I} = 0$, then $\chi\big(\T^k\backslash S_I^{2n + d - 1}\big) = 2d$.
\end{lemma}
\begin{proof}
If $I = \emptyset$, then $S_I^{2n + d - 1}$ is empty when $d = 0$ and consists of the two points $(0,\ldots,0,\pm 1)$, both fixed by $\T^k$, when $d = 1$,
So assume $I\neq\emptyset$ and $r_I = \abs{I} \geq 1$. As in the proof of Lemma~\ref{lem:RepRank=|I|},
we permute coordinates and row-reduce $A_I$ over $\Z$ to assume that $I = [r_I]$ and $A_I$ is upper-triangular, hence with nonzero entries on the diagonal.
First assume $d = 0$ and fix a point $(v_1,\ldots,v_n)\in S_I^{2n + d - 1}$. Proceeding as in the proof of Lemma~\ref{lem:RepRank=|I|}, we choose $t_1\in\T^1$
such that $t_1^{a_{1,1}} v_1$ is on the positive real axis, then choose $t_2\in\T^1$ such that $t_2^{a_{2,2}} t_1^{a_{1,1}} v_2$ is on the positive real
axis, etc., resulting in a $\mathbf{t}\in\T^k$ such that the coordinates of $\rho_I(\mathbf{t})(v_1,\ldots,v_n)$ are real, with the first $r_I$ coordinates positive and
the last $n - r_I$ coordinates zero. It follows that each point in $S_I^{2n + d - 1}$ is in the orbit of a point in the set
$\Delta_{r_I} = S_I^{2n + d - 1} \cap \R_+^{r_I} \times \{\mathbf{0}\}$ where $\R_+$ denotes the positive real axis.
Note that the map $(x_1,\ldots,x_{r_I},0,\ldots,0)\mapsto(x_1^2,\ldots,x_{r_I}^2)$
is a definable homeomorphism of $\Delta_{r_I}$ onto the interior of the standard $(r_I-1)$-simplex so that $\chi(\Delta_{r_I}) = (-1)^{r_I - 1}$.
As no two points with nonnegative real coordinates are in the same $\T^k$-orbit, we have that the quotient map $S_I^{2n + d - 1}\to\T^k\backslash S_I^{2n + d - 1}$
restricts to a continuous definable bijection $\Delta_{r_I}\to \T^k\backslash S_I^{2n + d - 1}$.
Then $\chi\big(\T^k\backslash S_I^{2n + d - 1}\big) = (-1)^{r_I - 1}$ follows from \cite[Ch.~4, Prop.~2.4]{vandenDriesBook},
the preservation of $\chi$ under definably bijections.

If $d = 1$, then for a point $(v_1,\ldots, v_n, w) \in S_I^{2n + d - 1}$, we apply the same procedure as above, obtaining a
$\mathbf{t}\in\T^k$ such that the first $r_I$ coordinates of $\rho_I(\mathbf{t})(v_1,\ldots,v_n)$ are positive real, the remaining coordinates zero.
Letting $S_{I,+}^{2n + d - 1}$ denote the set of points in $S_I^{2n + d - 1}$ with positive $w$-coordinate, it follows that
each orbit in $S_{I,+}^{2n + d - 1}$ contains a unique point with positive real $v_i$-coordinate for $i \leq r_I$ and positive $w$ coordinate.
As above, this space is homeomorphic to the standard $r_I$-simplex so that
$\chi\big(\T^k\backslash S_{I,+}^{2n + d - 1}\big) = (-1)^{r_I}$. 
The same argument shows that, for the space $S_{I,-}^{2n + d - 1}$ of points in $S_I^{2n + d - 1}$ with negative $w$-coordinate,
$\chi\big(\T^k\backslash S_{I,-}^{2n + d - 1}\big) = (-1)^{r_I}$, and for the space
$S_{I,0}^{2n + d - 1}$ of points in $S_I^{2n + d - 1}$ with $w$-coordinate zero,
$\chi\big(\T^k\backslash S_{I,0}^{2n + d - 1}\big) = (-1)^{r_I - 1}$. Hence
\[
    \chi\big(\T^k\backslash S_I^{2n + d - 1}\big) = (-1)^{r_I - 1} + 2(-1)^{r_I} = (-1)^{r_I}.
    \qedhere
\]
\end{proof}

We can now compute the $\Gamma$-Euler characteristic of $\T^k\ltimes S^{2n + d - 1}$.

\begin{theorem}
\label{thrm:SphereGamEC}
Let $A\in\Z^{k\times n}$ be the weight matrix of a unitary $\T^k$-representation $V$ and
let $U$ be a finite-dimensional real $\T^k$-representation with $U\simeq V$ or $U\simeq V\oplus W$ as above, where $\dim_\R W = 1$.
Let $\Gamma$ be a finitely presented group.
Choose a finite presentation of $\Gamma$ and let $N$ denote the number of generators and $r_M$ and $\sigma_1,\ldots,\sigma_{r_M}$ denote
the rank and invariant factors, respectively, of the matrix $M$ defined in Theorem~\ref{thrm:SmithAbelianization}.
Letting $d = \dim_\R U \bmod 2$ and $S^{2n + d - 1}$ denote the unit sphere in $U$,
\begin{equation}
\label{eq:SphereGamECUnitary}
\begin{split}
    \chi_\Gamma\big(\T^k\ltimes S^{2n + d - 1}\big)
    =
    \begin{cases}
        2d\prod\limits_{i=1}^{r_M} \sigma_i^k
            + \sum\limits_{\substack{\emptyset\neq I\subseteq[n] \\ r_I = \abs{I}}} (-1)^{r_I + d - 1}
                \prod\limits_{i=1}^{r_M} \sigma_i^{k - \abs{I}} \prod\limits_{j=1}^{\abs{I}} \gcd(\sigma_i, s_{I,j}),
            &   r_M = N,
        \\
        \sum\limits_{\substack{\emptyset\neq I\subseteq[n] \\ \abs{I} = k = r_I}} (-1)^{r_I + d - 1} \abs{\det A_I}^{N - r_M}
            \prod\limits_{i=1}^{r_M} \sigma_i^{k - \abs{I}} \prod\limits_{j=1}^{\abs{I}} \gcd(\sigma_i, s_{I,j}),
        &
        r_M \neq N,
    \end{cases}
\end{split}
\end{equation}
where $r_I = \Rank A_I$ and $s_{I,1}, \ldots, s_{I,r_I}$ are the invariant factors of $A_I$.
\end{theorem}
\begin{proof}
Applying Equation~\eqref{eq:GamECSumAbel} to the partition of $S^{2n + d - 1}$ into $\{ S_I^{2n + d - 1} : I\subseteq[n] \}$,
\[
    \chi_\Gamma\big(\T^k\ltimes S^{2n + d - 1}\big)
    =
    \sum\limits_{I\subseteq[n]} \chi\big(\T^k\backslash S_I^{2n + d -1}\big)\chi\big(\HOM(\Gamma, G_I)\big).
\]
The terms corresponding to $I$ such that $r_I < \abs{I}$ vanish by Lemma~\ref{lem:SphereRank<|I|}, so that by Lemma~\ref{lem:SphereRank=|I|},
\[
    \chi_\Gamma\big(\T^k\ltimes S^{2n + d - 1}\big)
    =
    2d\chi\big(\HOM(\Gamma, \T^k)\big)
        + \sum\limits_{\substack{\emptyset\neq I\subseteq[n] \\ r_I = \abs{I}}}
            (-1)^{r_I + d - 1}\chi\big(\HOM(\Gamma, G_I)\big).
\]
An application of Lemma~\ref{lem:RepGenGamHom} and the simplifications noted in Remark~\ref{rem:RepGamECSimplifications} complete the proof.
\end{proof}

% xxxxxxxxxxxxxxxxxxxxxxxxxxxxxxxxxxxxxxxxxxxxxxxxxxxxxxxxxxxxxxxxxxxxxxxxx

\subsection{Symplectic Quotients}
\label{subsec:GamECSubspacesSymplectic}

We continue to let $V$ be a unitary $\T^k$-representation with weight matrix $A$.
In this section, we consider the corresponding symplectic quotient associated to $V$.
Let us briefly describe the construction of the symplectic quotient; we refer the reader to \cite[Ch.~22-24]{daSilvaLecturesSympGeom},
\cite[Sec.~4.4]{GuilleminGinzburgKarshon}, and \cite{SjamaarLerman} for more details and note that we rescale the symplectic form as in
\cite[Sec.~2]{HerbigIyengarPflaum}. As a symplectic vector space admits a compatible complex structure, there is no need to consider a real
factor $W$ in this section.

Identifying the dual $(\mathfrak{t}^k)^\ast$ of the Lie algebra $\mathfrak{t}^k$ with $\R^k$,
the moment map $J\colon V\to(\mathfrak{t}^k)^\ast$ for the $\T^k$-action is given by
\begin{equation}
\label{eq:SympMoment}
    J(v_1,\ldots,v_n)   =   - A \begin{pmatrix} \abs{v_1}^2 \\ \vdots \\ \abs{v_n}^2 \end{pmatrix}
    \in(\mathfrak{t}^k)^\ast\simeq\R^k.
\end{equation}
Let $\mathbf{c} = (c_1,\ldots,c_k)\in\R^k$ denote an element of $J(V)\subset (\mathfrak{t}^k)^\ast$
and let $Z^{\mathbf{c}}$ denote $J^{-1}(\mathbf{c}) \subseteq V$. Note that $J$ is $\T^k$-invariant so that
$\T^k$ acts on $J^{-1}(\mathbf{c})$.
The \emph{symplectic quotient at level $\mathbf{c}$} is the quotient $\T^k\backslash Z^{\mathbf{c}}$.
We represent the symplectic quotient with the translation groupoid $\T^k\ltimes Z^{\mathbf{c}}$ and consider the $\Gamma$-Euler characteristic
$\chi_\Gamma\big(\T^k\ltimes Z^{\mathbf{c}}\big)$.
In \cite[Sec.~6]{FarsiMeitSeaton}, $\chi_\Gamma\big(G\ltimes Z^{\mathbf{0}}\big)$ was computed for an arbitrary compact Lie group $G$
and shown to be independent of the representation $V$; here, we are interested in the case $\mathbf{c}\neq\mathbf{0}$.

As in Section~\ref{sec:GamECFulLRep}, we will consider the decomposition of $Z^{\mathbf{c}}$ in terms of the subsets $I\subseteq[n]$.
Let $J_I$ denote the restriction of the moment map to $V_I$ and let $Z_I^{\mathbf{c}} = Z^{\mathbf{c}} \cap V_I$.
Clearly, we need only consider those $I$ such that $J_I = \mathbf{c}$ is consistent, i.e., such that $Z_I^{\mathbf{c}}\neq\emptyset$.
In particular, if $\mathbf{c} = \mathbf{0}$, then $Z_{\emptyset}^{\mathbf{c}}\neq\emptyset$.

Using an argument similar to that of Lemma~\ref{lem:RepRank<|I|}, we have the following.

\begin{lemma}
\label{lem:SymplecticRank<|I|}
Let $A\in\Z^{k\times n}$ be the weight matrix of a unitary $\T^k$-representation $V$.
Let $\mathbf{c}\in J(V)\subset (\mathfrak{t}^k)^\ast$ with $\mathbf{c}\neq\mathbf{0}$, and let $I\subseteq[n]$ such that $Z_I^{\mathbf{c}}\neq\emptyset$.
If $r_I < \abs{I}$, then $\chi\big(\T^k\backslash Z_I^{\mathbf{c}}\big) = 0$.
\end{lemma}
\begin{proof}
Assume $r_I < \abs{I}$ and $Z_I^{\mathbf{c}}\neq\emptyset$.
Let $T\simeq\T^{\abs{I}}$ denote the maximal compact torus of $\GL(V_I^{\cl},\C)$ as in the proofs of Lemmas~\ref{lem:RepRank<|I|} and \ref{lem:SphereRank<|I|}.
Once again, the action of $T$ on $V_I$ preserves the modulus of each coordinate so that as $J$ depends only on the moduli, $T$ acts on $Z_I^{\mathbf{c}}$.
As the coordinates of points in $Z_I^{\mathbf{c}}$ corresponding to $i\in I$ are nonzero, $(Z_I^{\mathbf{c}})^T = \emptyset$.
Then as $\rho_I(\T^k)$ is a subgroup of $T$ isomorphic to a torus of dimension $r_I < \abs{I}$, the claim follows from Proposition~\ref{prop:TorusActionSubgroupEC}.
\end{proof}

We also have an analog of Lemma~\ref{lem:RepRank=|I|} that is again proven similarly.

\begin{lemma}
\label{lem:SymplecticRank=|I|}
Let $A\in\Z^{k\times n}$ be the weight matrix of a unitary $\T^k$-representation $V$.
Let $\mathbf{c}\in J(V)\subset (\mathfrak{t}^k)^\ast$ with $\mathbf{c}\neq\mathbf{0}$, and let $I\subseteq[n]$ such that $Z_I^{\mathbf{c}}\neq\emptyset$.
If $r_I = \abs{I}$, then $\chi\big(\T^k\backslash Z_I^{\mathbf{c}}\big) = 1$.
\end{lemma}
\begin{proof}
Assume $r_I = \abs{I}$ and $Z_I^{\mathbf{c}}\neq\emptyset$.
As in the proofs of Lemma~\ref{lem:RepRank=|I|} and \ref{lem:SphereRank=|I|}, we permute coordinates and transform $\begin{pmatrix} A_I & \mathbf{c} \end{pmatrix}$
via an element of $\GL(k,\Z)$, corresponding to a change of coordinates for $\mathfrak{t}^k$ and $(\mathfrak{t}^k)^\ast$ and hence preserving $Z_I^{\mathbf{c}}$,
to assume that $I = [r_I]$ and $A_I$ upper-triangular, hence the diagonal entries are nonzero.
For a point $(v_1,\ldots,v_n)\in Z_I^{\mathbf{c}}$, it follows that the equation $J_I = \mathbf{c}$ uniquely determines $\abs{v_i}$ for
$i = 1,\ldots,r_I$, and $v_i = 0$ for $i > r_I$. Hence $Z_I^{\mathbf{c}}$ is homeomorphic to an $r_I$-dimensional real torus.
We once again choose, as in the proof of Lemma~\ref{lem:RepRank=|I|}, a $t_1\in\T^1$ such that $t_1^{a_{1,1}} v_1$ is on the positive real axis,
then $t_2\in\T^1$ such that $t_2^{a_{2,2}} t_1^{a_{1,1}} v_2$ is on the positive real axis, etc., resulting in a $\mathbf{t}\in\T^k$ such that
$\rho_I(\mathbf{t})(v_1,\ldots,v_n)$ is the unique point in $Z_I^{\mathbf{c}}$ whose first $r_I$ coordinates are positive real. It follows that each point in
$Z_I^{\mathbf{c}}$ is in the orbit of this point so that $\T^k\backslash Z_I^{\mathbf{c}}$ is a single point.
\end{proof}

We can now prove the following.

\begin{theorem}
\label{thrm:SymplecticGamEC}
Let $A\in\Z^{k\times n}$ be the weight matrix of a unitary $\T^k$-representation $V$ and let $\Gamma$ be a finitely presented group.
Choose a finite presentation of $\Gamma$ and let $N$ denote the number of generators and $r_M$ and $\sigma_1,\ldots,\sigma_{r_M}$ denote
the rank and invariant factors, respectively, of the matrix $M$ defined in Theorem~\ref{thrm:SmithAbelianization}. Then for each
$\mathbf{c}\in J(V)\subset (\mathfrak{t}^k)^\ast$ with $\mathbf{c}\neq\mathbf{0}$,
\begin{equation}
\label{eq:SymplecticGamEC}
    \chi_\Gamma\big(\T^k\ltimes Z^{\mathbf{c}}\big)
    =
    \sum\limits_{\substack{I\subseteq[n] \\ r_I = \abs{I}, \; Z_I^{\mathbf{c}}\neq\emptyset}} \abs{\det A_I}^{N - r_M}
        \prod\limits_{i=1}^{r_M} \sigma_i^{k - \abs{I}} \prod\limits_{j=1}^{\abs{I}} \gcd(\sigma_i, s_{I,j}),
\end{equation}
where $r_I = \Rank A_I$ and $s_{I,1}, \ldots, s_{I,r_I}$ are the invariant factors of $A_I$.
\end{theorem}
\begin{proof}
Applying Equation~\eqref{eq:GamECSumAbel} to the partition of $Z^{\mathbf{c}}$ into the $Z_I^{\mathbf{c}}$, we have
\[
    \chi_\Gamma\big(\T^k\ltimes Z^{\mathbf{c}}\big)
    =
    \sum\limits_{I\subseteq[n]} \chi\big(\T^k\backslash Z_I^{\mathbf{c}}\big)\chi\big(\HOM(\Gamma, G_I)\big).
\]
The terms corresponding to $I$ such that $Z_I^{\mathbf{c}} = \emptyset$ clearly vanish, as do the terms corresponding to $I$ such that
$r_I < \abs{I}$ by Lemma~\ref{lem:SymplecticRank<|I|}. This observation and Lemma~\ref{lem:SymplecticRank=|I|} yield
\[
    \chi_\Gamma\big(\T^k\ltimes Z^{\mathbf{c}}\big)
    =
    \sum\limits_{\substack{I\subseteq[n] \\ r_I = \abs{I}, \; Z_I^{\mathbf{c}}\neq\emptyset}}
        \chi\big(\HOM(\Gamma, G_I)\big).
\]
Then the result follows from Lemma~\ref{lem:RepGenGamHom}.
\end{proof}

It is interesting to observe that, unlike the symplectic quotient at level $\mathbf{0}$, the $\Gamma$-Euler characteristic of the
symplectic quotient at level $\mathbf{c}\neq\mathbf{0}$ depends on the representation and not merely on $\Gamma$ and $\T^k$.

\begin{example}
\label{ex:SymplecticNotRepIndy}
If $A = \diag(a_{1,1},\ldots,a_{n,n})$ is diagonal with each $a_{i,i} < 0$, then $J(\mathbf{v}) = \mathbf{c}$ is consistent for any
$\mathbf{c}\in(\mathfrak{t}^k)^\ast$ with positive coordinates, and by Theorem~\ref{thrm:SymplecticGamEC},
the $\Z^\ell$-Euler characteristic is given by $\abs{\det(A)}^\ell$ for any positive integer $\ell$.
In particular, for any $k \geq 1$, the representation can be chosen so that $\chi_\Z(\T^k\ltimes Z^{\mathbf{c}})$ has any positive integer value.
\end{example}

It is well known that if $\mathbf{c}$ is a regular value of the moment map, then the symplectic quotient $\T^k\backslash Z^{\mathbf{c}}$
is a symplectic manifold or orbifold \cite{MarsdenWeinsteinReduction,MeyerReduction}. In these cases, the $\Gamma$-Euler characteristic
computed by Theorem~\ref{thrm:SymplecticGamEC} coincides with the $\Gamma$-orbifold Euler characteristic defined in
\cite{TamanoiMorava,TamanoiCovering} for quotients by finite groups and \cite[Sec.~4]{FarsiSeaGenOrbEuler} for general orbifolds.
In many cases, however, $\T^k\backslash Z^{\mathbf{c}}$ has infinite isotropy groups and hence is not an orbifold, which we
observe with the following immediate consequence of Proposition~\ref{prop:SmithIsotropy}.

\begin{corollary}
\label{cor:SymplecticNotOrbifolds}
Let $A\in\Z^{k\times n}$ be the weight matrix of a unitary $\T^k$-representation $V$, and
let $\mathbf{c}\in (\mathfrak{t}^k)^\ast$. If there is an $I$ such that $r_I < k$ and $Z_I^{\mathbf{c}}\neq\emptyset$,
then points in $Z_I^{\mathbf{c}}$ have infinite isotropy. Hence, $\T^k\backslash Z^{\mathbf{c}}$ is not an orbifold.
\end{corollary}

% xxxxxxxxxxxxxxxxxxxxxxxxxxxxxxxxxxxxxxxxxxxxxxxxxxxxxxxxxxxxxxxxxxxxxxxxx
% xxxxxxxxxxxxxxxxxxxxxxxxxxxxxxxxxxxxxxxxxxxxxxxxxxxxxxxxxxxxxxxxxxxxxxxxx
% xxxxxxxxxxxxxxxxxxxxxxxxxxxxxxxxxxxxxxxxxxxxxxxxxxxxxxxxxxxxxxxxxxxxxxxxx

\section*{AI Disclosure}

The contents of this paper were human-generated. The second author used ChatGPT (OpenAI) and Claude (Anthropic) in the late stages of manuscript preparation
to identify references, fact-check and organize background information, and check completed drafts for correctness and typographical errors.
After using these tools, the authors addressed corrections manually and take full responsibility for the contents.

%------------------------------------------------------------------------------
%------------------------------------------------------------------------------
%------------------------------------------------------------------------------

\bibliographystyle{amsplain}
\bibliography{LRS-PS}

\end{document}